\documentclass[12pt]{article}

\usepackage{amsfonts,amsmath,amsthm,latexsym,amssymb,fullpage,amscd}
\usepackage[mathscr]{eucal}

\newcommand{\R}{\mathbb{R}}
\newtheorem{theorem}{Theorem}[section]
\newtheorem{lemma}[theorem]{Lemma}

\title{$L^p$ Sobolev regularity for a class of Radon and Radon-like transforms of various codimension}
\author{Michael Greenblatt}
\date{\today}

\newcommand\blfootnote[1]{%
  \begingroup
  \renewcommand\thefootnote{}\footnote{#1}%
  \addtocounter{footnote}{-1}%
  \endgroup
}
\begin{document}
\maketitle
\begin{abstract} 
In the paper [G1] the author proved  $L^p$ Sobolev regularity results for averaging operators over hypersurfaces and
connected them to associated Newton polyhedra. In this paper, we use rather different resolution of singularities techniques along with oscillatory
integral methods applied to surface measure Fourier transforms to
 prove $L^p$ Sobolev regularity results
for a class of averaging operators over surfaces which can be of any codimension.  
\end{abstract}

\section{ Background and the statement of Theorem 1.1.}

\blfootnote{ 2010 {\it Mathematics Subject Classification}: Primary: 42B20; Secondary: 44A12}
\blfootnote{{\it Key words}: Radon transform, Sobolev regularity, Oscillatory integral.}

Let $m \geq 1$ and $n > 1$ be integers. For $i= 1,...,n$ let $\gamma_i(t_1,...,t_m)$ be a nonconstant real analytic function on a neighborhood $U$ of the origin in $\R^m$ with $\gamma_i(0) = 0$ for each $i$. Assume that for each $i < n$, on the set $\{t \in U: \gamma_i(t) \neq 0\}$ one has
\[ \lim_{t \rightarrow 0} {\gamma_{i+1}(t) \over \gamma_i(t)} = 0 \tag{1.1} \]
We consider operators  of the following form, acting on functions on $\R^n$.
\[ {Tf(x) = \int_{\R^m} f(x  - \gamma(t))K(t)\,dt } \tag{1.2} \]
Here the kernel $K(t_1,...,t_m)$ is supported in $U$ and is $C^1$ on $U' = \{t \in U: t_i \neq 0$ for all $i\}$ such that for some 
$a_1,...,a_m$ with $0 \leq a_i < 1$ for each $i$ and some constant $C > 0$, on $U'$ one has the estimates
\[ |K(t_1,...,t_m)| \leq C|t_1|^{-a_1}...|t_m|^{-a_m} \tag{1.3a} \]
\[ |\partial_{t_i} K(t_1,...,t_m)| \leq C{1 \over |t_i|} |t_1|^{-a_1}...|t_m|^{-a_m} {\rm\,\,for\,\,all\,\, } i \tag{1.3b} \]
 
Examples of operators satisfying $(1.1)$ include when $\gamma_i(t) = \prod_{j=1}^i g_j(t)$ where $g_j(t)$ are nonconstant real analytic functions with $g_j(0) = 0$.
Another class of examples are $\gamma_i(t) = \sum_{j = 1}^m c_{ij} t_j^{2i}$ for positive constants $c_{ij}$. Also, whenever $m = 1$ one can do 
a linear transformation after which $(1.1)$ is satisfied since one can make the orders of the zeroes of the $\gamma_i(t)$ at $t = 0$ strictly increasing in $i$.
(Note that one might end up with one or more $\gamma_i(t)$ identically zero in this situation, which reduces things to a lower-dimensional 
scenario.)

The conditions $(1.3a)-(1.3b)$ on $K(t)$ describe $K(t)$ as a kind of multiparameter fractional integral kernel, and one
can often use resolution of singularities in the $t$ variables to 
reduce various other $K(t)$ to finite sums of operators where $(1.3a)-(1.3b)$ hold, as will be described at the end of this section. Traditional
Radon transforms (where $K(t)$ is smooth) fall under the situation where each $a_i = 0$ and our sharp estimates will also be sharp for such Radon 
transforms.

Our goal will be to prove $L^p(\R^n)$ to $L^p_{\beta}(\R^n)$ boundedness results  with $\beta > 0$ for the operator $(1.2)$ 
that are sharp up to endpoints for $p$ in a (not necessarily small) interval containing $p = 2$. Our theorem will be a local result, and will hold for $K(t)$ 
satisfying $(1.3a)-(1.3b)$ supported on a sufficiently small neighborhood of the origin. It follows from considerations similar to those of
 [CNStW] that for any 
$1 < p < \infty$ there will be some $\beta_p > 0$ such that $T$ is bounded from $L^p(\R^n)$ to $L^p_{\beta_p}(\R^n)$. This paper will describe up to endpoints,
for the class of averaging operators being considered here,  the optimal value of this $\beta_p$  for $p$ in an interval containing $p = 2$ and shows it is independent of $p$ in this interval, so long as this $\beta_p  < {1 \over n}$.

The level of smoothing in our theorems will be expressed in terms of the supremum $s_0$ of all $s$ such that there is a neighborhood $V$ of the 
origin such that $\int_{V}|\gamma_n(t)|^{-s}|t_1|^{-a_1} ... |t_m|^{-a_m}$ is finite. One can show using resolution of singularities that this supremum is 
positive and furthermore is independent of $V$ if $V$ is a small enough neighborhood of the origin. In fact one has a more precise statement. Let 
$d\mu$ denote the measure $|t_1|^{-a_1} ... |t_m|^{-a_m}\,dt$. Then there is an $s_0 > 0$ and an integer $0 \leq d_0 < m$ such that
 if $r$ is sufficiently small there is an asymptotic expansion of the following form as $\epsilon \rightarrow 0$, where $C_r \neq 0$.
\[\mu(\{t:|t| < r,\,|\gamma_n(t)| < \epsilon\}) = C_r \epsilon^{s_0} |\ln(\epsilon)|^{d_0} + o( \epsilon^{s_0} |\ln(\epsilon)|^{d_0}) \tag{1.4}\]
An elementary argument then shows that $s_0$ is  the supremum of the $s$ for which the integral
$\int_{|t| < r}|\gamma_n(t)|^{-s}\,|t_1|^{-a_1} ... |t_m|^{-a_m}\,dt$ is finite. We refer to the reference [AGuV] for more background on
 statements such as $(1.4)$ and related matters. 
 
 \noindent Our main theorem is as follows.\\
 
\begin{theorem} Let $T$ be as in $(1.2)$.

 \noindent {\bf 1)} Let $A$ be the open triangle with vertices $(0,0), (1,0)$, and $({1 \over 2}, {1 \over n})$, and let $B = \{(x,y) \in A: y < s_0\}$.
There is a neighborhood $V$ of the origin such that if $K(t)$ supported in $V$, then $T$ is bounded from $L^p(\R^n)$ to $L^p_{\beta}(\R^n)$ for $({1 \over p}, \beta)$ in $B$. 

\noindent {\bf 2)} Suppose $s_0 < 1$, $K(t)$ is nonnegative, and there exists a constant $C_1 > 0$ and a neighborhood $W$ of the origin such that $K(t) > C_1t_1^{-a_1}...t_m^{-a_m}$ for all  $t \in W$ with $t_i \neq 0$ for all $i$. Then if  $1 < p < \infty$ and $T$ is bounded from $L^p(\R^n)$ to
 $L^p_{\beta}(\R^n)$, we must have  $\beta \leq s_0$.
 
 \end{theorem}
 
 Observe that if $s_0 < {1 \over n}$, when taken together the two parts of Theorem 1.1 imply that for ${1 \over p} 
\in ( {ns_0 \over 2}, 1 -  {ns_0 \over 2})$, the amount of $L^p$ 
Sobolev smoothing given by part 1, $s_0$ derivatives, is optimal except possibly missing the endpoint $\beta = s_0$. When $s_0 = {1 \over n}$ the same is true for $p = 2$. A natural question to ask is whether or not the endpoint estimate $\beta = s_0$ also holds. 
It can be shown that when $p = 2$  the endpoint estimate
 holds if and only if $d_0 = 0$ in $(1.4)$. The author does not know what happens in the $p \neq 2$ situation.
 
 In a sense Theorem 1.1 will be proven by viewing a general operator $T$ satisfying $(1.1)$ and $(1.3a)-(1.3b)$ as an average of Radon transforms along curves (the $m=1$ case), using resolution of singularities to disentangle the curves appropriately. One then uses the Van der Corput lemma in one 
dimension on each of these curves in the blown up coordinates. Integrating the resulting estimates in the remaining $m-1$ 
variables will then give a desired result for
the original operator $T$. To obtain a result for $p \neq 2$, this process will be done in the context of an analytic interpolation.
 
As for the condition that $s_0 < 1$ in part 2) of Theorem 1.1, note that by $(1.1)$ the order of the zero of $\gamma_{i+1}(t)$ at the origin is
greater than the order of the zero of $\gamma_i(t)$ at the origin. So in particular $\gamma_n(t)$ has a zero of order at least $n$ at the origin, 
so that $|\gamma_n(t)| \leq C|t|^n$ for some constant $C$. Suppose
$m < n$, which is the usual situation when dealing with Radon transforms and their singular variants. Then $s_0$ is maximized when each $a_i = 0$ in 
$(1.3a)-(1.3b)$, in which case $s_0$ will be no greater than the value of $s_0$ if $\gamma_n(t)$ is replaced by $|t|^n$ in $(1.4)$. But a simple calculation
reveals $s_0 = {m \over n} < 1$ in this situation. Thus whenever $m < n$, in the setting of Theorem 1.1 we have $s_0 < 1$ and thus part 2) of 
Theorem 1.1 gives that the maximum possible amount of Sobolev smoothing for $T$ is given by $s_0$. 

\noindent {\bf Example 1.}

Consider the case where $m = 1$, so that $\gamma(t)$ is of the form $(c_1 t^{b_1}
+ o(t^{b_1}),..., c_n t^{b_n} + o(t^{b_n}))$ for some $b_1 < ... < b_n$ and some and nonzero $c_i$. Naturally, this example can be analyzed directly
without resolution of singularities, but this example is still useful in understanding the statement of Theorem 1.1.

The conditions $(1.3a)-(1.3b)$ translate into  $|K(t)| \leq Ct^{-a}$ and $\big|{dK \over dt} (t)\big| \leq Ct^{-a-1}$ for some 
$0 \leq a < 1$.
In view of $(1.4)$, $s_0$ is given in terms of the statement that if $r$ is sufficiently small, one has that $\mu(\{t: |t| < r, \,|t^{b_n}| < \epsilon\}) = 
c_r \epsilon^{s_0} |\ln (\epsilon)|^{d_0} + o(\epsilon^{s_0} |\ln (\epsilon)|^{d_0})$ for some $c_r \neq 0$. Since $d\mu = t^{-a}\,dt$, 
$s_0$ is given by the exponent in $\int_0^{\epsilon^{1 / b_n}} t^{-a} \,dt$, namely $s_0 = {1 - a \over b_n}$. So for example when $b_i = i$ for each $i$,
if $a = 0$ then Theorem 1.1 gives (other than the endpoint) the well known ${1 \over n}$ derivatives of $L^2$ Sobolev improvement one obtains for such
Radon transforms, while if $a > 0$ one obtains ${1 - a \over n}$ derivatives of Sobolev space improvement for a range of $p$ containing $p = 2$ that
grows as $a$ increases. 

It is worth pointing out that by [OSmSo], when $a = 0$ and $\gamma(t) = (\cos t, \sin t, t)$, which after a linear coordinate change falls into the  $b_1 = 1, b_2 = 2, b_3 = 3$ situation here, one obtains
Sobolev improvement in a larger range than the closure of the region provided by Theorem 1.1. 

\noindent {\bf Example 2.} 

Suppose each $a_i = 0$ and each $\gamma_i(t)$ is a monomial $t_1^{\alpha_{i1}}...t_m^{\alpha_{im}}$ such that the
multiindices $\alpha_i = (\alpha_{i1},...,\alpha_{im})$ are strictly increasing in $i$ under the lexicographical ordering. Then $s_0$ is given by the
condition that for small but fixed $r$ the Lebesgue measure of $\{t: |t| < r,\, |t^{\alpha_n}| \leq \epsilon\} \sim \epsilon^{s_0}|\ln(\epsilon)|^{d_0}$ as $\epsilon \rightarrow 0$. A direct calculation reveals that $s_0 = {1 \over \max_j \alpha_{nj}}$ and $d_0$ is one less than the number of times that the number $\alpha_{nj}$ appears in the multiindex $\alpha_n$.\\

\noindent {\bf Example 3.} 

Suppose $\gamma_i(t) = \sum_{j = 1}^m c_{ij} t_j^{2i}$ for positive constants $c_{ij}$. Then $s_0$ is determined by 
the condition that for small but fixed $r$ we have $\mu(\{t:|t| < r, \,\sum_{j = 1}^m c_{nj} t_j^{2n} < \epsilon\}) \sim \epsilon^{s_0}|\ln(\epsilon)|^{d_0}$ as $\epsilon \rightarrow 0$. Thus
$s_0$ is determined by $\int_{\{t:  |t_j| < \epsilon^{1 \over 2n}\forall j\}}\prod_{j=1}^m|t_j|^{-a_j}\,dt \sim \epsilon^{s_0}|\ln(\epsilon)|^{d_0}$.
Thus we have $s_0 = {m \over 2n} - {1 \over 2n} \sum_{j=1}^m a_j$.\\

\noindent {\bf Extensions to more general $K(t)$.} \\

Suppose that instead of satisfying $(1.3a)-(1.3b)$ the function $K(t)$ is of the form 
$\prod_{i=1}^M |g_i(t)|^{-\alpha_i}\beta(t)$ for a cutoff function $\beta(t)$ supported near the origin and some real analytic functions $g_i(t)$ on a neighborhood of the origin with $g_i(0) = 0$. Here $\alpha_i \geq 0$ such that the resulting $K(t)$ is integrable near the origin. Then one can reduce to the case of $K(t)$ satisfying $(1.3a)-(1.3b)$ using resolution of
singularities. For this, can use various forms of such theorems, including either Hironaka's famous work [H1] [H2] or the author's paper [G6]. If $\phi(t)$ is a smooth nonnegative bump function supported sufficiently near the origin and equal to one on a neighborhood of the 
origin, then by these resolution of singularities theorems one can write $\phi =\sum_{i=1}^N \phi_i(t)$ such that for each $i$ there is a coordinate change $\Psi_i$ with
$\Psi_i(0) = 0$ such that 
$\phi_i(\Psi_i(t))$ is smooth and supported on a neighborhood of the origin, and such that each $g_i(\Psi_i(x))$ and the Jacobian determinant of each $\Psi_i(t)$ are
both of the form $b_i(t)m_i(t)$ where $b_i(t)$ is smooth and nonvanishing on the support of $\phi_i ( \Psi_i(t))$, and where $m_i(t)$ is a monomial. 

Hence if $\beta(t)$ is supported in a neighborhood
 $V$ of the origin small enough that
$\phi(t) = 1$ on $V$, in $(1.2)$ one can write $K(t) = \sum_{i=1}^N K(t)\phi_i(t)$, and then do the coordinate changes $\Psi_i$ on the corresponding
integrals $(1.2)$. The result is that $Tf = \sum_{i=1}^N T_if$, where 
\[ T_i f(x)  = \int f(x - \gamma (\Psi_i(t))) \prod_{i=1}^M |g_i (\Psi_i (t))|^{-\alpha_i}|Jac(\Psi_i)(t)|\, \phi_i (\Psi_i(t)) \beta(\Psi_i(t))\,dt \tag{1.5}\]
Here $Jac(\Psi_i)(t)$ denotes the Jacobian determinant of $\Psi_i(t)$. Because $Jac(\Psi_i)(t)$ and the $g_i (\Psi_i(x))$ are monomialized in the above sense, the operator $T_i$ is of the form $(1.2)$ with a kernel satisfying $(1.3a)-(1.3b)$. The condition $(1.1)$ automatically holds in the blown up coordinates.

Suppose instead of satisfying $(1.3a)-(1.3b)$, $K(t)$ satisfies multiparameter fractional integral inequalities. In other words, suppose as in
the author's previous paper [G1]  we write 
$t = ({\bf t}_1,...,{\bf t}_p)$,
where ${\bf t}_i$ denotes  $(t_{i1},...,t_{il_i})$ such that the various $t_{ij}$ variables comprise the whole list $t_1,...,t_m$. Suppose there exist
$0 \leq \alpha_i < l_i$ and a constant $C > 0$ such that the following estimates hold:
\[|K(t)| \leq C\prod_{k=1}^p |{\bf t}_k|^{-\alpha_k} \tag{1.6a} \]
\[|\partial_{t_{ij}}K({\bf t})| \leq C {1 \over |t_{ij}|} \prod_{k=1}^p |{\bf t}_k|^{-\alpha_k} {\hskip 0.4 in} {\rm \,\,for\,\,all\,\,} i {\rm\,\, and\,\,} j 
\tag{1.6b}\]
Then one can do the above resolution of singularities process on the functions $|{\bf t}_k|^2$ to again reduce consideration to a finite sum of operators 
satisfying $(1.3a)-(1.3b)$. 

\noindent {\bf Some history.} 

There has been quite a bit of work on the boundedness properties of Radon transforms and related operators on various function spaces, so we focus
our attention on the $L^p$ Sobolev regularity questions being considered here. The following history is largely taken from [G1]. 

For translation-invariant Radon transforms over curves $\gamma(t) = (t, t^m)$ in the plane, $m > 2$, optimal $L^p$ to $L^p_{\beta}$ boundedness
properties are proven in [Gr] and [C]. These papers go beyond the results of this paper in that they also prove endpoint results. For example, when the function
$K(t)$ is smooth the results of [Gr] imply
that one has $L^p$ to $L^p_{\beta}$ boundedness in the closure of the region $A$ of Theorem 1.1 other than at $(1/p, \beta) = ({1 \over m}, {1 \over m})$
and $(1/p, \beta) = (1 - {1 \over m}, {1 \over m})$, and [C] shows that one does not have $L^p$ to $L^p_{\beta}$ boundedness for these two points.

For general non translation-invariant averages over curves in $\R^2$, again with smooth densities, thorough
$L^p_{\alpha}$ to $L^q_{\beta}$ boundedness results are proven in [Se] that are sharp up to endpoints. As mentioned in Example 1, in the case where $\gamma(t) = 
(\cos t, \sin t, t)$, [OSmSo] provides $L^p$ Sobolev regularity results that go beyond those that are provided by Theorem 1.1 after a linear change in coordinates.
 
For the case of translation-invariant Radon transforms over two-dimensional hypersurfaces in $\R^3$, there have been several papers. When $p = 2$,
a level of $L^2$ Sobolev improvement is equivalent to a statement on the decay rate of the Fourier transform of the surface measure in question, and
for the case of smooth density functions the stability results of Karpushkin [K1] [K2] provide sharp decay estimates. These estimates can then be interpreted
in terms of Newton polygons using [V]. Generalizations from real analytic to smooth surfaces then follow from [IKeM]. When one has a singular
density function, the paper [G5] proves some results that go beyond those that can be derived from those of this paper. Other results for singular
density functions appear in [G2].

For higher-dimensional hypersurfaces, in the situation of multiparameter fractional integral kernel $K(t)$, if the density functions are sufficiently singular then
there is an interval $I$ containing $p = 2$ such that sharp $L^p$ Sobolev improvement results for $p$ in $I$ follow from [S]. These results 
extend the results in [G4]. There is also the author's aforementioned paper [G1] which proves $L^p$ Sobolev regularity theorems for hypersurfaces when
$K(t)$ is a multiparameter fractional integral kernel. The paper [G1] also uses resolution of singularities but in quite a different way from this paper.
In the papers [G1] and [G5], when one gets estimates that are sharp up to endpoints analogously to Theorem 1.1, it can be shown that the amount
of Sobolev improvement given is exactly the quantity $s_0$ in Theorem 1.1, if one describes the hypersurface using the parameterization $\gamma(t) =  (t_1,...,t_{n-1}, \gamma_n(t_1,...,t_{n-1}))$ with $\gamma_n(0) = 0$ and $\nabla \gamma_n(0) = 0$. For the paper [G5] this is immediate from the definitions there, while for [G1] this can be shown in relatively short order using the results
of [G7].

Additional results for higher-dimensional hypersurfaces are in [Cu], and the paper [PSe] proves Sobolev space regularity theorems for certain
translation-invariant Radon transforms over curves. More generally, there have been several papers connecting $L^p$ Sobolev regularity of Radon transforms
to Fourier integral operators, including [GreSeW] [GreSe1] [GreSe2].

\section {The proof of part 1 of Theorem 1.1.}

\subsection {\bf Overview of argument.}

We will prove part 1 of Theorem 1.1 by embedding $T$ in an analytic family of operators as follows. We define the operator $T^z$ by 
\[T^z f(x)= e^{z^2}\int_{\R^m}|\gamma_n(t)|^z f(x  - \gamma(t))K(t)\,dt \tag{2.1}\]

We will readily see that if $0 > t_0 > \max(-s_0, -{1 \over n})$, then for $z$ on the line $Re(z) = t_0$ if  $2 < p < \infty$ one has estimates $||T^z f||_{L^p}
\leq C_{t_0}||f||_{L^p}$. We will then show that if $t_1 > \max(0,{1 \over n} - s_0)$, then on the line  $Re(z) = t_1$ one has estimates  $||T^z f||_{L_{{1 \over n}}^2}
\leq C_{t_1}||f||_{L^2}$. Interpolating between these two estimates gives a Sobolev space estimate for $T = T^0$.
 Then letting $p$ go to infinity, $t_0$ go to $ \max(-s_0, -{1 \over n})$, and $t_1$ go to $\max(0,{1 \over n} - s_0)$, part 1 of Theorem 1.1 will follow.
 
The argument for the line $Re(z) = t_0$ is straightforward. Since $T^z f$ is of the form  $f \ast \sigma_z$ for a surface-supported measure 
$\sigma_z$, we have that $||T^z||_{L^p \rightarrow L^p}$ is bounded by the integral of the magnitude of the density function in $(2.1)$, which in turn is bounded by
$C_{t_0}\int_{\R^m} |\gamma_n(t)|^{t_0}|K(t)|\,dt$. Since $s_0$ is defined to be the supremum of the $s$ for which the integral $ \int_{V} |\gamma_n(t)|^{-s}t_1^{-a_1}...t_m^{-a_m} \,dt$ is finite for sufficiently small neighborhoods of the origin, since $K(t)$ satisfies $(1.3a)$ and $t_0  > -s_0$  we have
that $C_{t_0}\int_{\R^m} |\gamma_n(t)|^{t_0}|K(t)|\,dt$ is finite. Hence $||T^z||_{L^p \rightarrow L^p}$ is uniformly bounded on the line
$Re(z) = t_0$ for any $p < \infty$ as needed.

\subsection  {$L^2$ estimates.}

Thus the main effort will go into proving the $L^2$ Sobolev estimates on the line
$Re(z) = t_1$. For our purposes it will not matter exactly which resolution of singularities procedure we use for the $L^2$ estimates and either
Hironaka's famous work [H1] [H2] or the author's [G6] will suffice. By these resolution of singularities results as applied simultaneously to the 
functions $\gamma_1(t),...,\gamma_n(t)$, there is an $r_0 > 0$ and
variable change maps $\Psi_j(t)$ such that  the following hold. Suppose $0 < r_2 < r_1 < r_0$ and $\phi(t)$ is a $C_c^\infty$ function such that 
$\phi(t)$ is supported on $\{t: |t| < r_1\}$ and such that $\phi(t) = 1$ on  $\{t: |t| < r_2\}$. Then $\phi(t)$ can be written as 
$\sum_{j=1}^M \phi_j(t)$ such that each $\phi_j (\Psi_j(t))$ is a  $C_c^\infty$ function on whose support each function
$\gamma_i (\Psi_j(t))$ can be written in the form $a_{ij}(t)m_{ij}(t)$, where $a_{ij}(t)$ doesn't vanish on the support of $\phi_j(\Psi_j(t))$
and $m_{ij}(t)$ is a nonconstant monomial $t_1^{\alpha_{ij1}}... t_m^{\alpha_{ijm}}$. It can also be arranged that $0$ is in the support of each
$\phi_j(\Psi_j(t))$ and that $\Psi_j(0) = 0$.

 Write $\alpha_{ij}$ to denote the multiindex $(\alpha_{ij1},...,\alpha_{ijm})$. If we had resolved the singularities of each $\gamma_{i_1}(t) - 
\gamma_{i_2}(t)$ for $i_1 \neq i_2$  simultaneously with the $\gamma_i(t)$, we automatically have that in the resolved coordinates for a given $j$ some
permutation of the multiindices 
$\{\alpha_{ij}\}_{i=1}^n$ is lexicographically ordered.  

In fact, $(1.1)$ ensures that we may further asssume that the multiindices $\alpha_{ij}$ are strictly lexicographically increasing in $i$ for a given $j$. To see why this is so, suppose
for some $j$ this were not the case, and there is an $i$ such that $\alpha_{ij} \geq \alpha_{i+1,j}$. Then since $\gamma_i (\Psi_j(t)) \sim 
t_1^{\alpha_{ij1}}... t_m^{\alpha_{ijm}}$ and $\gamma_{i+1}( \Psi_j(t))\sim 
t_1^{\alpha_{i+1j1}}... t_m^{\alpha_{i+1jm}}$, there is a constant $C > 0$ such that 
$|\gamma_{i+1} (\Psi_j(t))| > C|\gamma_i ( \Psi_j(t))|$ for all $t$ in the set $A$ defined as the set of points in the support of $\phi_j(\Psi_j(t))$ for
which  $t_k \neq 0$ for each $k$. In view of $(1.1)$, one must have that $\Psi_j(A)$ does not have the origin as a limit point, contradicting 
that $\Psi_j(0) = 0$ and that $0$ is in the support of $\phi_j(\Psi_j(t))$. Thus we may 
asssume that the multiindices $\alpha_{ij}$ are strictly lexicographically ordered in $i$ for a given $j$.

Let $a_j(t)m_j(t)$ denote $(a_{1j}(t)m_{1j}(t) ,..., a_{nj}(t)m_{nj}(t))$. Since $\phi(t) = 1$ on $|t| < r_2$, if the support of $K(t)$ is contained in
the set on which $\phi(t) = 1$ then we may write $K(t) = \sum_{j=1}^M K_j(t)$, where $K_j(t) = K(t)\phi_j(t)$. 
We write $T^zf = \sum_j T_j^z$ in accordance with the decomposition $K(t) = \sum_{j=1}^M K_j(t)$, and then perform
the coordinate changes $\Psi_j$, obtaining
\[T_j^z f(x) = e^{z^2} \int |a_{nj}(t)m_{nj}(t)|^z f(x - a_j(t)m_j(t)) Jac_j(t) K_j( \Psi_j(t)) \,dt \tag{2.2}\]
Here $Jac_j(t)$ denotes the Jacobian determinant of the coordinate change $\Psi_j(t)$. The resolution of singularities can always be done in such a way that $Jac_j(t)$ is also comparable to a monomial. 
If we had resolved each coordinate function $t_i$ simultaneously with the other functions in the above resolution of singularities process, then by the
chain rule each 
$ K_j( \Psi_j(t)) = K(t_1(\Psi_j(t)),...,t_n(\Psi_j(t)))\phi_j(\Psi_j(t))$ will 
satisfy $(1.3a)-(1.3b)$ as well, possibly with different $a_i$. Hence the same will be true of $L_j(t) = Jac_j(t) \times K_j ( \Psi_j(t))$. 
We now write $(2.2)$ succinctly as
\[T_j^z f(x) = e^{z^2} \int |a_{nj}(t)m_{nj}(t)|^z f(x - a_j(t)m_j(t)) L_j(t)\,dt \tag{2.3}\]
Here $L_j(t)$ satisfies $(1.3a)-(1.3b)$, possibly with different $a_i$. Our goal is to show for $Re(z) = t_1$ an estimate of the form
 $||T_j^z f||_{L_{{1 \over n}}^2} \leq C_{t_1}||f||_{L^2}$ for each $j$. Since $T_j^z f = f \ast \sigma_j^z$ for a measure $\sigma_j^z$, this in turn is
  equivalent to showing an estimate of
 the form $|\widehat{\sigma_j^z}(\lambda)| \leq C(1 + |\lambda|)^{-{1 \over n}}$. Writing out $\widehat{\sigma_j^z}(\lambda)$, we have
\[\widehat{\sigma_j^z}(\lambda) = e^{z^2} \int |a_{nj}(t)m_{nj}(t)|^z e^{-i\lambda \cdot a_j(t)m_j(t)} L_j(t)\,dt \tag{2.4}\]
 It will simplify our arguments if for any given $j$ there is a single $k$ for which the exponent $\alpha_{ijk}$ is strictly increasing in $k$. This 
 will automatically be the case for $m = 1$. In order to ensure this is the case when $m > 1$, we perform an additional (relatively simple) resolution of singularities. We proceed as follows. We divide $\R^m$
 (up to a set of measure zero) into $m!$ regions $\{A_k\}_{k=1}^{m!}$, where 
each $A_k$ is a region of the form $\{t \in \R^m: |t_{l_1}| < |t_{l_2}| < ... < |t_{l_m}|\}$, where  $t_{l_1},...,t_{l_m}$ is a permutation of
the $t$ variables. We write $(2.4)$ as $\sum_{k=1}^{m!} I_{jk}^z(\lambda)$,
where
\[I_{jk}^z(\lambda) = e^{z^2} \int_{A_k} |a_{nj}(t)m_{nj}(t)|^z e^{-i\lambda \cdot a_j(t)m_j(t)} L_j(t)\,dt \tag{2.5}\]
We focus our attention one one such $A_k$ and we let $u_i$ denote $t_{l_i}$. We make the variable changes $u_i = \prod_{p=i}^m y_p$. In the $y$ variables, $A_k$  becomes the rectangular box $(-1,1)^{n-1} \times \R$ and $(2.5)$ can be written in the form
\[I_{jk}^z(\lambda) = e^{z^2} \int_{(-1,1)^{n-1} \times \R} |b_{njk}(y)p_{njk}(y)|^z e^{-i\lambda \cdot b_{jk}(y)p_{jk}(y)} M_{jk}(y)\,dy \tag{2.6}\]
Here $b_{jk}(y)p_{jk}(y)$ is of the form $(b_{1jk}(y)p_{1jk}(y),..., b_{njk}(y)p_{njk}(y))$, where 
the $p_{ijk}(y)$ are again monomials, the $b_{ijk}(y)$ are nonvanishing real-analytic functions, and $M_{jk}(y)$ again satisfy $(1.3a)-(1.3b)$, but with 
different exponents. Note that each $p_{ijk}(y)$ is of the form $y_1^{\beta_{ijk1}}...y_m^{\beta_{ijkm}}$ where $\beta_{ijkm}$ is the overall degree of the monomial 
$m_{ij}(t)$. Since the multiindices $\alpha_{ij}$ are strictly increasing, the exponent $\beta_{ijkm}$ is strictly increasing in $i$ for fixed $j$ and $k$.

The above considerations were for the situation when $m > 1$, but we can incorporate the $m = 1$ case into the above by simply letting there be one
$k$ for each $j$ and take the $y$ variable to be the original $t$ variable. So in the following we assume no restrictions on $m$ and $(2.6)$ will still 
hold in this sense.

The idea behind the analysis of $(2.6)$ is as follows. For fixed values of $y_1,...,y_{m-1}$, the phase function $\lambda \cdot b_{jk}(y)p_{jk}(y)$ is 
effectively of the form $c_1y_m^{\beta_{1jkm}} + ... + c_ny_m^{\beta_{njkm}}$ for $c_i = \lambda_i b_{ijk}(y_1,...,y_{k-1},0)y_1^{\beta_{ijk1}}...
y_{m-1}^{\beta_{ijk\,m-1}}$ and thus the Van der Corput lemma can be 
used in the $y_m$ variable since the exponents $\beta_{ijkm}$ are distinct for $i = 1,...,n$. 

To simplify our notation, we let $\rho_i$ denote $\beta_{ijkm}$.  Write $g(y_m) = \sum_{i=1}^n
 c_i y_m^{\rho_i}$. Since the $\rho_i$ are distinct, the vectors $\{(\rho_1^p,...,\rho_n^p): p = 1,...,n\}$ are linearly independent (as 
can be shown using the Vandermonde determinant). As a result, using row operations the vectors $w_p = (\prod_{q=0}^{p-1}(\rho_1 - q), ...,
\prod_{q=0}^{p-1} (\rho_n - q) )$ are linearly independent for $p = 1,...,n$. Consequently, there exists an $\epsilon > 0$ such that for
any vector $v$ in $\R^n$ there is some $p$ with $1 \leq p \leq n$ such that
\[|w_p \cdot v| \geq \epsilon|v| \tag{2.7}\]
Letting $v = c_iy_m^{\rho_i}$, equation $(2.7)$ implies that for each $y_m$ there necessarily exists a $p$ with $1 \leq p \leq n$ (which can depend on $y_m$) such that
\[\bigg|{d^p g \over d y_m^p}(y_m)\bigg| > \epsilon{1 \over |y_m|^p} \sum_{i=1}^n |c_i y_m^{\rho_i}| \tag{2.8}\]
Since there is a constant $C$ such that $\big|{d^{p+1} g \over d y_m^{p+1}}(y_m)\big| < C{1 \over |y_m|^{p+1}} \sum_{i=1}^n |c_i y_m^{\rho_i}|$
for all $p$ and all $y$, in view of $(2.8)$ each dyadic interval in $y_m$ can be written as the union of boundedly many subintervals on each of which we have
for some $1 \leq p \leq n$ that
\[\bigg|{d^p g \over d y_m^p}(y_m)\bigg| > {\epsilon \over 2} {1 \over |y_m|^p} \sum_{i=1}^n |c_i y_m^{\rho_i}| \tag{2.9}\]
Going back to $(2.6)$, in order to use the Van der Corput Lemma we would like $(2.9)$ to hold with $g(t)$ replaced by the phase function in $(2.6)$. That is,
we would like $(2.9)$ to hold if one replaces $c_i = \lambda_i b_{ijk}(y_1,...,y_{m-1},0)y_1^{\beta_{ijk1}}... y_{m-1}^{\beta_{ijk\,m-1}}$ with
the function of $y_m$ given by $\lambda_i b_{ijk}(y_1,...,y_{m-1},y_m)y_1^{\beta_{ijk1}}... y_{m-1}^{\beta_{ijk\,m-1}}$. Since each $b_{ijk}(y_1,...,y_{m-1},y_m)$ is real analytic and nonvanishing, this will hold if $|M_{jk}(y)|$ is supported on $|y_m| < \delta$ for $\delta$ sufficiently small; 
the effect of a $y_m$ derivative landing on $b_{ijk}(y_1,...,y_{m-1},y_m)$ is to introduce a bounded factor, which will be much smaller than the  $C {1 \over |y_m|}$ factor that drives
the estimate $(2.9)$ if $|y_m|$ is sufficiently small.

So the question becomes whether or not we can assume $|y_m| < \delta$ for $\delta$ small enough so that the analogue of $(2.9)$ holds for the phase function in $(2.6)$. To see why the answer is yes, first note that in terms of the $t$ variables, one will have that each $|y_m| < \delta$ if one has  
$|t| < \delta'$ for some $\delta'$ depending on $\delta$. 

Next, let $t'$ be in the support of some $\phi_j(\Psi_j(t))$ such that $\Psi_j(t') = 0$. There is a neighborhood of $t'$ on which we may shift 
coordinates to become centered at $t = t'$ instead of
$t = 0$; the functions that were monomialized before will be 
monomialized in the shifted coordinates, and the multiindices in the shifted coordinates will be lexicographically ordered in the shifted coordinates for the
same reason as before. Note also that this shift may change one or more $a_p$ in $(1.3a)-(1.3b)$ into zero if the $p$th component of $t_{jl}'$ is nonzero.

We let $D_{t'}$ be a disk centered at $t'$ such that  the above considerations hold on $D_{t'}$ and such that 
each $|t|$ is less than the associated $\delta'$ in the coordinates centered at $t'$. By compactness, we
may let $W_j$ be a neighborhood of the points in the support of $\phi_j(\Psi_j(t))$ with $\Psi_j(t) = 0$ that is a finite union of such $D_{t'}$. There is 
then an $r_{W_j} > 0$
such that if $\bar{\phi}(t)$ is a smooth function supported on $|t| < r_{W_j}$ and equal to 1 on a neighborhood of the origin, then if we write $\bar{\phi}_j(t) = \bar{\phi}(t)\phi_j(t)$, then $\sum_j \bar{\phi}_j(t) =  \bar{\phi}(t)$ and $\bar{\phi}_j(\Psi_j(t))$ is a smooth function supported in $W_j$. 

We can then use a partition of unity subordinate to the $D_{t'}$ comprising $W_j$ to write 
$\bar{\phi}_j(\Psi_j(t))$ as a finite sum $\sum_l h_{jl}(t)$ of smooth functions such that each $h_{jl}(t)$ is supported on one of the $D_{t'}$. We denote
the $t'$ corresponding to $h_{jl}(t)$ by $t_{jl}'$.
If we use the $h_{jl}(t + t_{jl}')$ in place of ${\phi}_j(\Psi_j(t))$ and the $h_{jl} (\Psi^{-1}(t) - t_{jl}'))$ in place
of ${\phi}_j(t)$, then the associated $|y_m|$ will be as small as needed above. Thus if $\bar{\phi}(t)$ is supported on $\{t: |t| < \min_j r_{W_j}\}$, 
then if the neighborhood $V$ on which $K(t)$ is supported is
contained in the set on which $\bar{\phi}(t)$ = 1, one can write $K(t) = \sum_{jl}K_{jl}(t)$ where $K_{jl}(t) = K(t)h_{jl} (\Psi^{-1}(t) - t_{jl}')$
and decompose $T = \sum_{jl} T_{jl}$ accordingly. For each $j$ and $l$, in the blown up and shifted coordinates the $|y_m|$ will be smaller
 than $\delta'$ as desired. 
 
 Thus replacing our original decompositon $T = \sum_j T_j$ by the decomoposition $T = \sum_{jl} T_{jl}$ if necessary, we
 assume that we always have $|y_m| < \delta'$. As a result, we can assume that the $|y_m|$ are small enough 
for fixed $(y_1,...,y_{m-1})$ such that each
dyadic interval in the $y_m$ variable can be written as the union of boundedly many subintervals on which for some $1 \leq p \leq n$ we have
\[\bigg|{\partial_{y_m}^p (\lambda \cdot b_{jk}(y)p_{jk}(y)) \over \partial y_m^p}\bigg| > {\epsilon \over 4} {1 \over |y_m|^p} \sum_{i=1}^n |c_i y_m^{\rho_i}| \tag{2.10}\]
On each of these intervals in the $y_m$ variable, we will apply the Van der Corput lemma in conjunction with $(2.10)$ in the integral $(2.6)$. 
When $p > 1$, we use the standard Van der Corput lemma (see p. 334 of [St]):\\

\begin{lemma}
Suppose $h(x)$ is a $C^k$ real-valued function on the interval $[a,b]$ such that $|h^{(k)}(x)| > A$ on $[a,b]$ for
some $A > 0$. Let $\phi(x)$ be $C^1$ on $[a,b]$. If $k \geq 2$ there is a constant $c_k$ depending only on $k$ such that
\[\bigg|\int_a^b e^{ih(x)}\phi(x)\,dx\bigg| \leq c_kA^{-{1 \over k}}\bigg(|\phi(b)| + \int_a^b |\phi'(x)|\,dx\bigg) \tag{2.11}\]
If $k =1$, the same is true if we also assume that $h(x)$ is $C^2$ and $h'(x)$ is monotone on $[a,b]$. 
\end{lemma}

\noindent  When $p = 1$, we will make use of the following variant of Lemma 2.1 for $k = 1$.\\

\begin{lemma} (Lemma 2.2 of [G3].) Suppose the hypotheses of Lemma 2.1 hold with $k = 1$, except  instead of assuming that $h'(x)$ is monotone on $[a,b]$ we assume that $|h''(x)|
< {B  \over(b-a)}A$ for some constant $B > 0$. Then we have
\[\bigg|\int_a^b e^{ih(x)}\phi(x)\,dx\bigg| \leq  A^{-1}\bigg((B+2) \sup_{[a,b]}|\phi(x)| + \int_a^b |\phi'(x)|\,dx \bigg) \tag{2.12}\]
\end{lemma}

\noindent In the integral $(2.6)$ in the $y_m$ variable, the function denoted by $\phi(x)$ in $(2.11)$ and $(2.12)$ is given by $|b_{njk}(y)p_{njk}(y)|^z  M_{jk}(y)$.
Recall that $M_{jk}(y)$ satisfies $(1.3a)-(1.3b)$ except with different exponents. So for some $b_1,...,b_m$ we have
\[|M_{jk}(y)| \leq C|y_1|^{-b_1}...|y_m|^{-b_m} \tag{2.13a}\]
\[ |\partial_{y_i} M_{jk}(y)| \leq C{1 \over |y_i|} |y_1|^{-b_1}...|y_m|^{-b_m} {\rm\,\,for\,\,all\,\, } i \tag{2.13b}\]
Recall also that $p_{njk}(y)$ is a monomial and that $b_{njk}(y)$ is
real analytic function which doesn't vanish on the support of the integrand of $(2.6)$. Hence for some constant $C' > 0$ we have
\[\big| |b_{njk}(y)p_{njk}(y)|^z  M_{jk}(y) \big| \leq C' |p_{njk}(y)|^{{\rm\,Re}(z)}|y_1|^{-b_1}...|y_m|^{-b_m} \tag{2.14a}\]
\[\big|\partial_{y_m}  |b_{njk}(y)p_{njk}(y)|^z  M_{jk}(y) \big| \leq C' |z| {1 \over |y_m|} |p_{njk}(y)|^{{\rm\,Re}(z)}|y_1|^{-b_1}...|y_m|^{-b_m} \tag{2.14b}\]
For fixed $y_1,...,y_{m-1}$, we now apply Lemma 2.1 or 2.2 in the $y_m$ directions on one of the boundedly many subintervals of some dyadic 
interval $J' = 2^{-q-1} \leq |y_m| < 2^{-q}$ on which $(2.9)$ holds. Denote this subinterval by $J$ and let $y_m'$ denote the center of this interval.
Then we obtain
\[\bigg|\int_J |b_{njk}(y)p_{njk}(y)|^z e^{i\lambda \cdot b_{jk}(y)p_{jk}(y)} M_{jk}(y)\,dy_m\bigg| \leq \]
\[C''|z|
|p_{njk}(y_1,...,y_m')|^{{\rm\,Re}(z)} \times |y_m'| \times |y_1|^{-b_1}...|y_m'|^{-b_m} \times (\sum_{i=1}^n |c_i (y_m')^{\rho_i}|)^{-{1 \over p}} \tag{2.15a}\]
Since $y_m \sim y_m'$ on $J'$, this can be reexpressed as 
\[\bigg|\int_J |b_{njk}(y)p_{njk}(y)|^z e^{i\lambda \cdot b_{jk}(y)p_{jk}(y)} M_{jk}(y)\,dy_m\bigg| \leq \]
\[C'''|z|\int_{J'}
|p_{njk}(y_1,...,y_m)|^{{\rm\,Re}(z)}\times |y_1|^{-b_1}...|y_m|^{-b_m}\times (\sum_{i=1}^n |c_i y_m^{\rho_i}|)^{-{1 \over p}}\,dy_m\tag{2.15b}\]
Notice that on $J'$, $c_i = \lambda_i b_{ijk}(y_1,...,y_{m-1},0)y_1^{\beta_{ijk1}}... y_{m-1}^{\beta_{ijk\,m-1}}$ is comparable in magnitude to 
$\lambda_i b_{ijk}(y_1,...,y_{m-1},y_m)y_1^{\beta_{ijk1}}... y_{m-1}^{\beta_{ijk\,m-1}}$. Since $\rho_i$ was defined to be $\beta_{ijkm}$, the term 
$|c_i y_m^{\rho_i}|$ is comparable in magnitude to $\lambda_i b_{ijk}(y_1,...,y_{m-1},y_m)y_1^{\beta_{ijk1}}... y_m^{\beta_{ijkm}}$. Note
that there is some $i$ for which $|\lambda_i| > {1 \over n}|\lambda|$, and since the multiindices $(\beta_{ijk1},...,\beta_{ijkm})$ are lexicographically
increasing in $i$, this $|\lambda_i b_{ijk}(y_1,...,y_m)y_1^{\beta_{ijk1}}... y_m^{\beta_{ijkm}}|$ is bounded below by a constant times
$|\lambda| | b_{njk}(y_1,...,y_m)y_1^{\beta_{njk1}}... y_m^{\beta_{njkm}}|$. Since $ b_{njk}(y_1,...,y_m)$ is nonvanishing on the support of
the integrand in question, we also have a lower bound of a constant times $|\lambda| |y_1^{\beta_{njk1}}... y_m^{\beta_{njkm}}|$. Hence we have
\[\sum_{i = 1}^n |c_i y_m^{\rho_i}| \geq C''''|\lambda| |y_1^{\beta_{njk1}}... y_m^{\beta_{njkm}}| \tag{2.16}\]
As a result, the right-hand side of $(2.15b)$ is bounded by
\[C_0|z| \int_{J'}
|p_{njk}(y)|^{{\rm\,Re}(z)}|y_1|^{-b_1}...|y_m|^{-b_m}(|\lambda| |y_1^{\beta_{njk1}}... y_m^{\beta_{njkm}}|)^{-{1 \over p}}\,dy_m\tag{2.17}\]
Simply by taking absolute values of the integrand and integrating, we also have
\[\bigg|\int_J |b_{njk}(y)p_{njk}(y)|^z e^{i\lambda \cdot b_{jk}(y)p_{jk}(y)} M_{jk}(y)\,dy_{m}\bigg| 
\leq C_1\int_{J' } |p_{njk}(y)|^{{\rm\,Re}(z)}|y_1|^{-b_1}...|y_m|^{-b_m}\,dy_m \tag{2.18}\]
Combining $(2.17)$ and $(2.18)$ we get a bound of
\[\leq C_2(1 + |z|)\int_{J'} |p_{njk}(y)|^{{\rm\,Re}(z)}|y_1|^{-b_1}...|y_m|^{-b_m}
\min(1,(|\lambda| |y_1^{\beta_{njk1}}... y_m^{\beta_{njkm}}|)^{-{1 \over p}})\,dy \tag{2.19}\]
Since $1 \leq p \leq n$, the right-hand side of $(2.19)$ is maximized for $p = n$. Inserting $p = n$ and adding $(2.19)$ over the boundedly many 
subintervals $J$ corresponding to a given $J'$ we obtain
\[\bigg|\int_{J'}|b_{njk}(y)p_{njk}(y)|^z e^{i\lambda \cdot b_{jk}(y)p_{jk}(y)} M_{jk}(y)\,dy_{m}\bigg|\]
\[\leq C_2(1 + |z|)\int_{J'} |p_{njk}(y)|^{{\rm\,Re}(z)}|y_1|^{-b_1}...|y_m|^{-b_m} \min(1,(|\lambda| |y_1^{\beta_{njk1}}...
 y_m^{\beta_{njkm}}|)^{-{1 \over p}})\,dy \tag{2.20}\]
Adding $(2.20)$ over all $J'$ for which the integrand is not identically zero, and then integrating the result in 
the $y_1,...,y_{m-1}$ variables leads to the following for some $\delta_0 > 0$.
\[\bigg|\int_{(-1,1)^{n-1} \times \R} |b_{njk}(y)p_{njk}(y)|^z e^{i\lambda \cdot b_{jk}(y)p_{jk}(y)} M_{jk}(y)\,dy\bigg|\]
\[\leq C_2(1 + |z|)\int_{(-1,1)^{n-1} \times [-\delta_0, \delta_0] } |p_{njk}(y)|^{{\rm\,Re}(z)}|y_1|^{-b_1}...|y_m|^{-b_m}\min(1, (|\lambda| |y_1^{\beta_{njk1}}... y_m^{\beta_{njkm}}|)^{-{1 \over p}})\,dy \tag{2.21}\]
As a result, for some $\delta_0 > 0$ we have the following bound on $I_{jk}^z$ of equation $(2.6)$.
\[|I_{jk}^z| \leq C_3 (1 + |z|) e^{Re(z)^2 - Im(z)^2} \int_{(-1,1)^{n-1} \times [-\delta_0,\delta_0] } |p_{njk}(y)|^{{\rm\,Re}(z)}|y_1|^{-b_1}...|y_m|^{-b_m}\]
\[\times \min(1,(|\lambda| |y_1^{\beta_{njk1}}... y_m^{\beta_{njkm}}|)^{-{1 \over n}})\,dy \tag{2.22}\]
Note that $(1 + |z|) e^{Re(z)^2 - Im(z)^2} $ is uniformly bounded in $Im(z)$ for fixed $Re(z)$. Hence we can replace $C_3 (1 + |z|) 
e^{Re(z)^2 - Im(z)^2}$ by $C_{t_1}$ where $t_1$ denotes $Re(z)$. Furthermore, the form of $(2.22)$ is such that if we replace $\delta_0$ by a smaller
$\delta_1$, inequality $(2.22)$ will still hold, but with a different constant $C_{t_1,\delta_1}$. So we have
\[|I_{jk}^z| \leq C_{t_1,\delta_1} \int_{(-1,1)^{n-1} \times [-\delta_1,\delta_1] } |p_{njk}(y)|^{t_1}|y_1|^{-b_1}...|y_m|^{-b_m}\]
\[\times \min(1,(|\lambda| |y_1^{\beta_{njk1}}... y_m^{\beta_{njkm}}|)^{-{1 \over n}})\,dy \tag{2.23}\]
In particular, we may assume $\delta_1$ is small enough so that the pullback of $(-1,1)^{n-1} \times [-\delta_1,\delta_1]$ under the coordinate changes
of the above resolution of singularities is contained in a set $\{t: |t| < r\}$ on which $(1.4)$ is valid. Converting back into the original $t$ coordinates through
these coordinate changes, recalling that  $C_1 < |b_{njk}(y_1,...,y_m)| < C_2$  for some $C_1, C_2 > 0 $, equation $(2.23)$ then gives 
\[|I_{jk}^z| \leq C_{t_1}'\int_{\{t: |t| < r\}} |\gamma_n(t)|^{t_1}|t_1|^{-a_1}...|t_m|^{-a_m}\min(1,(|\lambda| |\gamma_n(t)|)^{-{1 \over n}})\,dy \tag{2.24a}\]
In the notation of $(1.4)$, $(2.24a)$ is simply  
\[|I_{jk}^z| \leq C_{t_1}'\int_{\{t: |t| < r\}} |\gamma_n(t)|^{t_1}\min(1,(|\lambda| |\gamma_n(t)|)^{-{1 \over n}})\,d\mu \tag{2.24b}\]
Recall that we are trying to show that $|I_{jk}^z| \leq C_{t_1}(1 + |\lambda|)^{-{1 \over n}}$ for $t_1 > max(0, {1 \over n} - s_0)$, where $s_0$ is as 
in $(1.4)$. If $|\lambda| < 1$ this is immediate from taking absolute values of the integrand in $(2.6)$ and integrating, so we assume $|\lambda| > 1$.
It is natural to break the integral in $(2.24b)$ into $|\gamma_n(t)| < {1 \over |\lambda|}$ and $|\gamma_n(t)| \geq {1 \over |\lambda|}$ parts. 
Equation $(2.24b)$ then becomes
\[|I_{jk}^z| \leq C_{t_1}'\int_{\{t: |t|< r,\, |\gamma_n(t)|  \leq {1 \over |\lambda|}\}}|\gamma_n(t)|^{t_1}\,d\mu + C_{t_1}'|\lambda|^{-{1 \over n}}
\int_{\{t: |t|< r,\,|\gamma_n(t)|  \leq {1 \over |\lambda|}\}}|\gamma_n(t)|^{t_1 - {1 \over n}}\,d\mu \tag{2.25}\]
Since $t_1 > {1 \over n} - s_0$, the exponent $t_1 - {1 \over n}$ in the second integral of $(2.25)$ is greater than $-s_0$. Since $s_0$ satisfies $(1.4)$,
this integral has a finite bound independent of $\lambda$. Hence the second term in $(2.25)$ has the desired bound of $C_{t_1}''|\lambda|^{-{1 \over n}}$.
As for the first term of $(2.25)$, we write it as
\[C_{t_1}' \sum_{j = 0}^{\infty} \int_{\{t:|t|< r,\, \leq {2^{-j - 1} \over |\lambda|} \leq |\gamma_n(t)|  \leq {2^{-j} \over |\lambda|}\}}|\gamma_n(t)|^{t_1}\,d\mu \tag{2.26}\]
In a given term of $(2.26)$, $|\gamma_n(t)|^{t_1} \leq |\lambda|^{-t_1}2^{-jt_1}$, and by $(1.4)$ for any $s < s_0$ the measure of the domain of integration of 
$(2.26)$ is bounded by $C_s|\lambda|^{-s}2^{-js}$. Thus $(2.26)$ is bounded by $C_{t_1,s} |\lambda|^{-s - t_1}\sum_{j=0}^{\infty} 2^{-j(s + t_1)}
= C_{t_1,s}'|\lambda|^{-s-t_1}$. Since $t_1 > {1 \over n} - s_0$, if $s$ is close enough to $s_0$ we have $s + t_1 > {1 \over n}$, and we have a bound 
of $C_{t_1}''|\lambda|^{-{1 \over n} - \epsilon}$ for some $\epsilon > 0$, better than what is needed.

Thus we conclude that $|I_{jk}^z| \leq C_{t_1}'''(1 + |\lambda|)^{-{1 \over n}}$ for some constant $C_{t_1}'''$ as desired. Adding this over all $j$ and $k$ 
then gives the estimate $|\widehat{\sigma_z}(\lambda)| \leq  C_{t_1}''''(1 + |\lambda|)^{-{1 \over n}}$, where $\sigma_z$ is such that $T^z f = 
f \ast \sigma_z$. Consequently we obtain the desired estimate $||T^z||_{L^2 \rightarrow L^2_{1 \over n}} \leq C_{t_1}'''''$ whenever $t_1 = 
Re(z) > \max(0, {1 \over n} - s_0)$.

\subsection {The end of the proof of part 1 of  Theorem 1.1.}

The argument here is very similar to the interpolation arguments in [G1] and [G5].
We just saw that if  ${\rm Re}\,(z) = t_1 > \max(0,{1\over n} - s_0)$ we have $||T^z f||_{L^2_{1 \over n}} \leq C_{t_1}||f||_{L^2}$.
At the beginning of this section, we saw that if $p < \infty$ and ${\rm\,Re}(z) = t_0 > \max(-s_0, -{1 \over n})$, then we have an estimate
$||T^z f||_{L^p} \leq C_{t_0}||f||_{L^p}$.

 Note that
$0 = \alpha\max(-s_0, -{1 \over n}) +  (1 - \alpha)(\max(0,{1\over n} - s_0)) $, where $\alpha = \max(0, 1 - s_0 n)$. 
Thus  if $1 > \alpha' > \alpha$,
 one can write $0 = \alpha' t_0 + (1 - \alpha') t_1$, where $0 > t_0 >\max(-s_0, -{1 \over n})$ and  $t_1  > \max(0,{1\over n} - s_0)$. Hence by  complex interpolation $T = T^0$ is bounded from $L^q$ to $L^q_{\beta}$, where
 ${1 \over q} =  \alpha' {1 \over p} + (1 -  \alpha') {1 \over 2}$ and $\beta = \alpha' 0 + (1 - \alpha'){1 \over n}$. Explicitly, we have
 $q = {1 \over {1 \over 2} + \alpha'({1 \over p} - {1 \over 2})}$ and $\beta = {1 - \alpha' \over n}$. 
 
 Using interpolation again, we have that $T$ is bounded from $L^r$ to $L^r_{\gamma}$ for $({1 \over r}, \gamma)$ in the closed 
 triangle with vertices $(0,0), (1,0),$ and $({1 \over q}, \beta)$. Taking the union of these triangles as $\alpha'$ approaches $\alpha$ and $p$ approaches
 $\infty$, we get that
 $T$ is bounded from $L^r$ to $L^r_{\gamma}$  for $({1 \over r}, \gamma)$ in the open triangle with vertices $(0,0), (1,0),$ and $({1 \over q'}, \beta')$, where 
$q' =  {1 \over {1 \over 2} - {1 \over 2}\alpha} = {2 \over \min(1,s_0n)} = \max(2, {2 \over s_0 n})$ and where $\beta' = 
 {1 - \alpha \over  n} = \min({1 \over n}, s_0)$.
 
 In the case where $s_0  \geq {1 \over  n}$, the union of these triangles is the open triangle with vertices $(0,0)$, $(1,0)$, and $({1 \over 2},
 {1 \over n})$, which is the boundedness region stipulated by part 1 of Theorem 1.1 in this case.
 If $s_0 < {1 \over n}$, 
 the union of these triangles  is the open triangle with vertices $(0,0)$, $(1,0)$, and $({s_0 n \over 2}, s_0)$. So $T$ is bounded from $L^r$ to $L^r_{\gamma}$ for $({1 \over r}, \gamma)$ in this region. By duality, it is also  bounded from $L^r$ to $L^r_{\gamma}$ for $({1 \over r}, \gamma)$ 
 such that $(1 - {1 \over r}, \gamma)$ is in this region, giving the triangle with vertices $(0,0)$, $(1,0)$, and $(1 - {s_0 n \over 2}, s_0)$. 
 Thus $T$ is bounded from $L^r$ to $L^r_{\gamma}$ for $({1 \over r}, \gamma)$ in the open trapezoidal region with vertices $(0,0)$, $(1,0)$, 
 $({s_0 n \over 2}, s_0)$, and $(1 - {s_0 n \over 2}, s_0)$. This is the region stipulated by Theorem 1.1 in the case where  $s_0 < {1 \over  n}$. This completes the proof of part 1 of Theorem 1.1.
 
 \section{The proof of part 2  of Theorem 1.1.}

 The proof of the second part of Theorem 1.1 is much easier than the first part and does not require resolution of singularities. Effectively
the argument reduces to showing (up to endpoints) the well-known fact that the sublevel set measure growth rate is at least as fast as 
the scalar oscillatory integral decay 
rate when the latter index is less than 1. While resolution of singularities can be used for this part, we instead use a more direct approach very similar 
to the corresponding argument in [G1].

 We suppose the hypotheses of part 2 of Theorem 1.1 hold. That is, we assume that $s_0 < 1$, $K(t)$ is nonnegative, and there exists a constant $C_1 > 0$ and a neighborhood $W$ of the origin such that $K(t) > C_1t_1^{-a_1}...t_m^{-a_m}$ for all $t \in W$ with $t_i \neq 0$ for all $i$. 
 Suppose that  $T$ is bounded from $L^p(\R^n)$ to $L^p_{\beta}(\R^n)$ for some $ 1 < p < \infty$. Then by duality, $T$ is bounded from $L^q(\R^n)$ to $L^q_{\beta}(\R^n)$ where $q$ is such that ${1 \over p} + {1 \over q} = 1$. 
 
 Since either $p \leq 2 \leq q$ or $q \leq 2 \leq p$, using interpolation
 we have that 
$T$ is bounded from $L^2(\R^n)$ to $L^2_{\beta}(\R^n)$. As a result, if $\sigma$ denotes the measure such that $Tf = f \ast \sigma$, we have the estimate
$|\hat{\sigma}(\lambda)| \leq C(1 + |\lambda|)^{-\beta}$. Explicitly, this means that
\[\bigg|\int_{\R^m} e^{-i\lambda \cdot \gamma(t)} K(t)\,dt \bigg| \leq C(1 + |\lambda|)^{-\beta}
\tag{3.1}\]
In particular, $(3.1)$ holds in the $(0,...,0,1)$ direction. So for all $\tau \in \R$ we have
\[\bigg|\int_{\R^m} e^{i\tau \gamma_n(t)} K(t)\,d t\bigg| \leq C(1 + |\tau|)^{-\beta} \tag{3.2}\]
Denote the integral on the left of $(3.2)$ by $U(\tau)$.
Let $B(x)$ be a bump function on $\R$  whose Fourier transform is nonnegative, compactly supported, and equal to 1 on a 
neighborhood of the origin, and let $\epsilon$ be a small  positive number. If $0 <  \beta' < \min(\beta,1)$, then $(3.2)$ implies
 that for some constant $A$ independent of $\epsilon$ one has
\[\int_{\R}|U(\tau) \tau^{\beta '  - 1}B(\epsilon \tau)|\,d \tau < A \tag{3.3}\]
Inserting the definition of $U(\tau)$ we have
\[\bigg|\int_{\R^{m+1}} e^{i\tau \gamma_n(t)} K(t) |\tau|^{\beta'  - 1}B(\epsilon \tau)\,d\tau \,dt\bigg|  <  A \tag{3.4}\]
We do the integral in $\tau$ in $(3.4)$. Letting $b_{\epsilon}(y)$ denote the convolution of $|y|^{-\beta '}$ with ${1 \over \epsilon} \hat{B}({y \over \epsilon})$,  for a constant $A'$ independent of $\epsilon$ we get
\[\bigg|\int_{\R^m}b_{\epsilon} (-\gamma_n(t)) \,K(t)\,dt\bigg| < A' \tag{3.5}\]
Note that both $b_{\epsilon} (-\gamma_n(t))$ and $K(t)$ are nonnegative here. Thus we may remove the absolute value and let
 $\epsilon \rightarrow 0$ to obtain
\[\int_{\R^m}|\gamma_n(t)|^{-\beta '} K(t) < \infty  \tag{3.6}\]
Since $ K(t)$ is bounded below by $ C_1t_1^{-a_1}...t_m^{-a_m}$ on a neighborhood $W$ of the origin, we therefore have
\[\int_{W}|\gamma_n(t)|^{-\beta '}\prod_{k=1}^m |t|^{-a_k} \,dt< \infty  \tag{3.7}\]
In other words, $|\gamma_n(t)|^{-\beta '}$ is in $L^1(W)$ with respect to the measure $\mu$. Hence it is in weak $L^1$, 
and we have the existence of a constant $C$ such that
	\[\mu (\{t \in W : |\gamma_n(t)|^{-\beta '} > \epsilon \}) \leq C  {1 \over  \epsilon} \tag{3.8}\]
Replacing $ \epsilon$ by $ \epsilon^{-\beta '}$, we get
\[\mu (\{t \in W: |\gamma_n(t)| <  \epsilon \}) \leq C  \epsilon^{\beta '} \tag{3.9}\]
In view of $(1.4)$, we have $\beta ' \leq s_0$. Since this holds for each $\beta'$ satisfying
$0 < \beta' <  \min(\beta, 1)$, we conclude that $\min(\beta,1) \leq s_0$. Since we are assuming $s_0 < 1$, we obtain that
$\beta \leq s_0$ as needed. This completes the proof of part 2 of Theorem 1.1.

\section {References.}

\noindent [AGuV] V. Arnold, S. Gusein-Zade, A. Varchenko, {\it Singularities of differentiable maps},
Volume II, Birkhauser, Basel, 1988. \setlength{\parskip}{0.3 em}

\noindent [C] M. Christ, {\it Failure of an endpoint estimate for integrals along curves} in Fourier analysis and partial differential equations (Miraflores de la Sierra, 1992), 163-168, Stud. Adv. Math., CRC, Boca Raton, FL, 1995.

\noindent [CNStW] M. Christ, A. Nagel, E. M. Stein, and S. Wainger, {\it Singular and maximal Radon transforms: analysis and geometry},
Ann. of Math. (2) {\bf 150} (1999), no. 2, 489-577. 

\noindent [Cu] S. Cuccagna, {\it  Sobolev estimates for fractional and singular Radon transforms}, J. Funct. Anal. {\bf 139} (1996), no. 1, 94-118.

\noindent [Gr] L. Grafakos, {\it Endpoint bounds for an analytic family of Hilbert transforms}, Duke Math. J. {\bf 62} (1991), no. 1, 23-59.

\noindent [G1] M. Greenblatt, {\it $L^p$ Sobolev regularity of averaging operators over hypersurfaces and the Newton polyhedron}, J. Funct. Anal. (2018)
https://doi.org/10.1016/j.jfa.2018.05.014

\noindent [G2] M. Greenblatt, {\it Uniform bounds for Fourier transforms of surface measures in $\R^3$ with nonsmooth density},
Trans. Amer. Math. Soc. {\bf 368} (2016), no. 9, 6601-6625.

\noindent [G3] M. Greenblatt, {\it Fourier transforms of irregular mixed homogeneous hypersurface measures},  Math. Nachr. 291 (2018), no. 7, 1075-1087. 

\noindent [G4] M. Greenblatt, {\it An analogue to a theorem of Fefferman and Phong for averaging operators
along curves with singular fractional integral kernel}, Geom. Funct. Anal. {\bf 17} (2007), no. 4, 1106-1138.

\noindent [G5] M. Greenblatt, {\it Smooth and singular maximal averages over 2D hypersurfaces and associated Radon transforms}, submitted.

\noindent [G6] M. Greenblatt, {\it A constructive elementary method for local resolution of singularities}, preprint.

\noindent [G7] M. Greenblatt, {\it Oscillatory integral decay, sublevel set growth, and the Newton
polyhedron}, Math. Annalen {\bf 346} (2010), no. 4, 857-895.

\noindent [GreSe1] A. Greenleaf, A. Seeger, {\it  Fourier integral operators with fold singularities}, J. Reine Angew. Math. {\bf 455} (1994), 35-56.

\noindent [GreSe2]  A. Greenleaf, A. Seeger, {\it Fourier integral operators with cusp singularities},  Amer. J. Math. {\bf 120} (1998), no. 5, 1077-1119.

\noindent [GreSeW] A. Greenleaf, A. Seeger, S. Wainger, {\it Estimates for generalized Radon transforms in three and four dimensions}, in Analysis, geometry, number theory: the mathematics of Leon Ehrenpreis (Philadelphia, PA, 1998), 243-254, Contemp. Math., {\bf 251}, Amer. Math. Soc., Providence, RI, 2000.

\noindent [H1] H. Hironaka, {\it Resolution of singularities of an algebraic variety over a field of characteristic zero I}, 
 Ann. of Math. (2) {\bf 79} (1964), 109-203.

\noindent [H2] H. Hironaka, {\it Resolution of singularities of an algebraic variety over a field of characteristic zero II},  
Ann. of Math. (2) {\bf 79} (1964), 205-326. 

\noindent [IKeM] I. Ikromov, M. Kempe, and D. M\"uller, {\it Estimates for maximal functions associated
to hypersurfaces in $\R^3$ and related problems of harmonic analysis}, Acta Math. {\bf 204} (2010), no. 2,
151-271.

\noindent [K1] V. N. Karpushkin, {\it A theorem concerning uniform estimates of oscillatory integrals when
the phase is a function of two variables}, J. Soviet Math. {\bf 35} (1986), 2809-2826.

\noindent [K2] V. N. Karpushkin, {\it Uniform estimates of oscillatory integrals with parabolic or 
hyperbolic phases}, J. Soviet Math. {\bf 33} (1986), 1159-1188.

\noindent [OSmSo] Daniel Oberlin, Hart Smith, and Christopher Sogge, {\it Averages over curves with torsion.} (English summary)
Math. Res. Lett. {\bf 5} (1998), no. 4, 535-539.

\noindent [PSe] M. Pramanik, A. Seeger, {\it $L^p$ Sobolev regularity of a restricted X-ray transform in $\R^3$} (English summary) Harmonic analysis and its applications, 47-64, Yokohama Publ., Yokohama, 2006. 

\noindent [Se] A. Seeger, {\it Radon transforms and finite type conditions}, J. Amer. Math. Soc. {\bf 11} (1998), no. 4, 869-897.

\noindent [St] E. Stein, {\it Harmonic analysis; real-variable methods, orthogonality, and oscillatory \hfill\break
integrals}, Princeton Mathematics Series Vol. 43, Princeton University Press, Princeton, NJ, 1993.

\noindent [S] B. Street, {\it Sobolev spaces associated to singular and fractional Radon transforms}, Rev. Mat. Iberoam. {\bf 33} (2017), no. 2, 633-748.

\noindent [V] A. N. Varchenko, {\it Newton polyhedra and estimates of oscillatory integrals}, Functional 
Anal. Appl. {\bf 18} (1976), no. 3, 175-196.
\\
\\

\noindent Department of Mathematics, Statistics, and Computer Science \hfill \break
\noindent University of Illinois at Chicago \hfill \break
\noindent 322 Science and Engineering Offices \hfill \break
\noindent 851 S. Morgan Street \hfill \break
\noindent Chicago, IL 60607-7045 \hfill \break
\noindent greenbla@uic.edu \hfill\break

\end{document}